\documentclass[twoside,11pt]{article}
\title{} \author{} \date{}
\usepackage{latexsym,amssymb,times}
\newtheorem{te}{Theorem}[section]

\newtheorem{cor}[te]{Corollary}
\newtheorem{fac}[te]{Fact}

\newtheorem{rem}[te]{Remark}

\def\dok{\noindent{\bf Proof. }}
\def\kdok{\hfill $\Box$ \par \vspace*{2mm} }

\def\f{\varphi}
\def\o{\omega}
\def\k{\kappa}
\def\r{\mathrel{\rho }}
\def\s{\mathrel{\sigma }}
\def\e{\varepsilon}

\def\G{{\mathbb G}}

\def\X{{\mathbb X}}

\def\N{{\mathbb N}}
\def\X{{\mathbb X}}

\def\CU{{\mathcal U}}

\def\la{\langle}
\def\ra{\rangle}

\begin{document}
\thispagestyle{plain}

\begin{center}
           {\large \bf CONNECTED REDUCED PRODUCTS}
\end{center}
\begin{center}
{\bf Milo\v s S.\ Kurili\'c}\footnote{Department of Mathematics and Informatics, Faculty of Sciences, University of Novi Sad,
                                      Trg Dositeja Obradovi\'ca 4, 21000 Novi Sad, Serbia.
                                      e-mail: milos@dmi.uns.ac.rs}
\end{center}
\begin{abstract}
\noindent
A structure $\X=\la X,\r\ra$, where $\r$ is a binary relation on the set $X$, is connected
iff the minimal equivalence relation containing $\r$ is the full relation on $X$.
Connectivity is not a first order property but it is expressible by a $L_{\o _1\o}$-sentence.
We show that, for a set $I$ the following conditions are equivalent
\begin{itemize}
\item[(a)] $|I|$ is less than the first measurable cardinal,
\item[(b)] For each filter $\Phi \subset P(I)$ and each family  $\{ \X_i :i\in I\}$ of  binary structures,
the reduced product $\prod _\Phi \X _i$ is connected,
iff there are a finite set $K\subset I$ and $n\in \o$ such that $\X _i$ is connected, for each $i\in K$
and  $\{ i\in I: \X _i \mbox{ is of diameter }\leq n\}\cup K\in \Phi$,
\item[(c)] The ultraproduct $\prod _{\CU}\G_\o$ is a disconnected graph for each non-principal ultrafilter $\CU \subset P(I)$, where $\G_\o$ is the linear graph on $\o$.
\end{itemize}
Moreover, the implication ``$\Leftarrow$" in (b) holds in ZFC.
For a disconnected $\prod _\Phi \X _i$ we describe the components of $\prod _\Phi \X _i$.


\end{abstract}
\section{Preliminaries}
\paragraph{Reduced products}
Let $L =\{ R \}$ be the language with one binary relational symbol.
Let $I$ be a non-empty set, $\Phi \subset P(I)$ a filter and $\X _i =\la X_i ,\r _i \ra$, for $i\in I$, $L$-structures.
Let $\sim _\Phi$ be the binary relation on the product $\prod _{i\in I}X _i$ defined in the following way:
for $x=\la x_i :i\in I\ra, \,y=\la y_i :i\in I\ra\in \prod _{i\in I}X _i$
$$
x \sim _\Phi y \Leftrightarrow \{ i\in I : x_i =y_i\}\in \Phi .
$$
It is evident that the relation $\sim _\Phi$ is an equivalence relation on the set $\prod _{i\in I}X _i$.
Let $[x]_\Phi =\{ y\in \prod _{i\in I}X _i :x \sim _\Phi y\}\in (\prod _{i\in I}X _i)/\sim _\Phi$ denote the equivalence class of $x$.
The corresponding {\it reduced product} $(\prod _{i\in I}\X _i )/\Phi$ (shortly $\prod _\Phi \X _i $) is the structure
$$\textstyle
\Big\la (\prod _{i\in I}X _i)/\!\sim _\Phi ,\r \Big\ra ,
$$
where the binary relation $\r$ on the quotient $(\prod _{i\in I}X _i)/\!\sim _\Phi$ is defined by
$$
[x]_\Phi \r  [y]_\Phi \Leftrightarrow \{ i\in I : x_i\r _i y_i\}\in \Phi .
$$
We will use the following evident fact (see \cite{Chang}, p.\ 177).
\begin{fac}\label{T014}
If $\Phi$ is a filter on $I$ and $J\in \Phi$, then $\Phi \upharpoonright J:=\Phi \cap P(J)$ is a filter on $J$ and
$f :\prod _\Phi \X _i \cong \prod _{\Phi \upharpoonright J} \X _j:=(\prod _{i\in J}\X _i)/\!\sim_{\Phi \upharpoonright J}$,
where $f ([x]_\Phi)=[x \upharpoonright J]_{\Phi \upharpoonright J}$.
\end{fac}
First-order properties of reduced products are well understood. If $L$ is any language,
an $L$-formula is a {\it basic Horn formula} iff it is of the form $\eta _0 \lor \dots \lor \eta _{m-1}$,
where $\eta _i$'s are literals and at most one of them is atomic.
The set of {\it Horn formulas} is the closure of the set of basic Horn formulas under $\land$, $\exists$ and $\forall$.
It is known that each  Horn formula $\f (v_0, \dots, v_{n-1})$ is {\it preserved in reduced products};
i.e., for each set $I$, filter $\Phi \subset P(I)$, family of $L$-structures $\{ \X _i :i\in I\}$
and $n$-tuple $\la x^0, \dots ,x^{n-1}\ra\in (\prod _{i\in I}X _i)^n$ we have
\begin{equation}\label{EQ028} \textstyle
\Big\{ i\in I : \X _i \models \f [x^0_i, \dots ,x^{n-1}_i]\Big\}\in \Phi \;\Rightarrow \;\prod _\Phi \X _i \models \f \Big[ [x^0]_\Phi, \dots ,[x^{n-1}]_\Phi \Big] .
\end{equation}
By the well known result of Keisler (under CH) and Galvin (in ZFC)  we have:
an $L$-sentence $\f$ is {\it preserved in reduced products}, i.e., for each $I$, $\Phi$
and $\{ \X _i :i\in I\}$ as above
\begin{equation} \label{EQ029}\textstyle
\{ i\in I : \X _i \models \f\}\in \Phi \;\Rightarrow \;\prod _\Phi \X _i \models \f ,
\end{equation}
iff $\f$ is logically equivalent to a Horn $L$-sentence (see \cite{Chang}, p.\ 365).

Putting the implication ``$\Leftarrow$" in (\ref{EQ029}) we obtain a definition of a sentence {\it preserved under reduced factors}.
It is known that the positive sentences (without $\neg$) are of that kind (see \cite{Chang}, p.\ 340).
\paragraph{Connected binary structures}
Let $L =\{ R \}$ be the language with one binary relational symbol.
Let $\X =\la X ,\r\ra$ be an $L$-structure, where $\rho \subset X^2$ is a reflexive binary relation  and let $\r ^{-1}$ be the inverse of $\r$. For $\e \in 2:=\{0,1\}$ let
\begin{equation}\label{EQ045}
\r ^\e=\left\{\begin{array}{ll}
                   \r       & \mbox{ if } \e = 0,\\
                   \r ^{-1} & \mbox{ if } \e = 1.
               \end{array}
         \right.
\end{equation}
If $x,y\in X$ and $n\in \o$, a {\it path of length $n+1$ from $x$ to $y$} is a pair $\la \bar z,\e\ra$,
where $\bar z=\la z_0,\dots,z_{n-1}\ra\in X^n$ and $\e =\la \e_0,\dots,\e_{n}\ra\in 2^{n+1}$ and such that
$$
x\r ^{\e _0}z_0\r ^{\e _1}z_1 \r ^{\e _2} \dots \r ^{\e _{n-1}}z_{n-1}\r ^{\e _n}y.
$$
(So, for $n=0$, $x\r ^{\e _0}y$.)
Then we will write $d(x,y)\leq n+1$
and $d(x,y)=\infty$ will denote that there is no (finite) path from $x$ to $y$.

The binary relation $\sim$ on the set $X$ defined by $x\sim y$ iff there is  finite path from $x$ to $y$,
is an equivalence relation,
the corresponding equivalence classes are the {\it (connected) components of $\X$}
and the structure $\X$ is {\it connected} iff it has only one component.
(It is easy to see that  $\sim$ is the transitive closure $(\r _s)_t$ of the symmetrization $\r _s:=\r \cup \r ^{-1}$ of the relation $\r$.)

We extend these definitions for arbitrary $L$-structure $\X =\la X ,\r\ra$ in the following way:
instead of $\rho $ we take its reflexivization $\rho _R:=\rho \cup \Delta _X$ and call $\X$ connected
iff  $\X _R :=\la X ,\rho _R\ra$ is connected,  components of $\X$ are components of $\X _R$, etc.

The statement $d(x,y)\leq  n+1$ is expressible by the first order formula
$$
\f _{d\leq n+1}(x,y) :=\exists z_0,\dots,z_{n-1}\bigvee _{\e \in 2^{n+1}}x\r ^{\e _0}z_0\r ^{\e _1} \dots \r ^{\e _{n-1}}z_{n-1}\r ^{\e _n}y.
$$
It is evident that $\f _{d\leq n+1}$ is positive and $\neg \f _{d\leq n+1}$ is Horn.
The property of an $L$-structure that each two points are connected by a path of length $\leq n+1$
is expressible by the first order sentence
$$
\f _n^{\rm conn}:=\forall x,y \;\;\f _{d\leq n+1}(x,y).
$$
If such $n\in \o$ exists we will say that $\X$ is {\it of finite diameter}.
The connectedness of an $L$-structure is expressible by the $L_{\o _1\o}$-sentence
$$
\f ^{\rm conn}:=\forall x,y \;\;\bigvee _{n\in \o }\f _{d\leq n+1}(x,y).
$$

\section{Connected reduced products. Components}
In the sequel we assume that
\begin{itemize}
\item[($\ast$)] $I\neq \emptyset$ is a set, $\Phi $ is a filter in $P(I)$, $\X _i =\la X_i, \r _i\ra$, $i\in I$, are $L$-structures
and $\prod \X _i :=\la \prod _{i\in I}X_i , \s \ra$ and $\prod _{\Phi }\X _i :=\la \prod _{i\in I}X_i / \!\!\sim _{\Phi }, \r\ra$
are the corresponding direct and reduced product respectively.
\end{itemize}
\begin{fac}\label{T017}
{\rm ($\ast$)}
If $x,y\in \prod _{i\in I}X_i$, $A\subset I$, $x\upharpoonright (I\setminus A)=y\upharpoonright (I\setminus A)$, $n\in \o$,
and there is $\e \in 2^{n+1}$ such that for each $i\in A$ there is $\bar z\in X^n_i$ such that $\la \bar z,\e\ra$
is a path from $x_i$ to $y_i$, then in $\prod _{\Phi }\X _i$ there is a path from $[x]_\Phi$ to $[y]_\Phi$.
\end{fac}
\dok
Let $\e=\la \e _0 ,\dots,\e _n\ra$, and for $i\in A$ let  $\la \la z^0_i, \dots,z^{n-1}_i \ra,\e \ra$ be a path from $x_i$ to $y_i$
\begin{equation}\label{EQ400}
x_i\r ^{\e_0}_{i} z^0_i\r ^{\e_1}_{i}z^1_i \r ^{\e_2}_{i} \dots \r ^{\e_{n-1}}_{i}z^{n-1}_i\r ^{\e_n}_{i}y_i.
\end{equation}
Let $t^m\in \prod _{i\in I}X_i$, for $m<n$, be defined by
\begin{equation}\label{EQ057}
t^m_i=\left\{\begin{array}{ll}
                   x_i & \mbox{ if } i\in I\setminus A ,\\
                   z^m_i& \mbox{ if } i\in A.
               \end{array}
         \right.
\end{equation}
By (\ref{EQ400}) in the direct product we have $x \s^{\e_0} t^0 \s^{\e_1} t^1\s^{\e_2} \dots \s^{\e_{n-1}} t^{n-1}\s ^{\e_n}y$
and, hence, in the reduced product $\prod _\Phi \X _i$  we have
\begin{equation}\label{EQ049}
[x]_\Phi\r ^{\e_0}  [t^0]_\Phi\r ^{\e_1} [t^1]_\Phi \r ^{\e_2}
\dots \r ^{\e_{n-1}} [t^{n-1}]_\Phi\r ^{\e_n} [y]_\Phi.
\end{equation}
So, this is a path from $[x]_\Phi$ to $[y]_\Phi$.
\kdok

\begin{fac}\label{T019}
If $\k\geq \o$ is a non-measurable cardinal and $\Phi \subset P(\k)$ a filter consisting of infinite sets,
then there is a partition $\k=\bigcup _{n\in \o}I_n$ such that
for each $F\in \Phi$ we have $F\cap I_n\neq \emptyset$, for infinitely many $n\in \o$.
\end{fac}
\dok
First let $\bigcap \Phi=\emptyset$.
Let $\CU \subset P(\k)$ be an ultrafilter such that $\Phi \subset \CU$.
Then the ultrafilter $\CU$ is non-principal
and, since the cardinal $\k$ is non-measurable, $\CU$ is not $\o$-complete.
Thus there is a partition $\k=\bigcup _{n\in \o}I_n$ such that $I_n\not\in \CU$
and, hence, $\k \setminus I_n \in \CU$, for each $n\in \o$.
Assuming that there are $F\in \Phi$ and $m\in \o$ such that $F\subset \bigcup _{n\leq m}I_n$,
we would have $F \cap \bigcap _{n\leq m}\k \setminus I_n=F \setminus \bigcup _{n\leq m} I_n=\emptyset \in \CU$,
which is false.
Thus, $|\{ n\in \o :F\cap I_n\neq \emptyset \}|=\o$, for each $F\in \Phi$.

Second let $A:=\bigcap \Phi\neq\emptyset$.
If $|A|< \o$, then, since $\Phi \subset P(\k)\setminus [\k ]^{<\o}$, $\Phi':=\Phi \upharpoonright (\k \setminus A)$
is a non-principal filter on $\k \setminus A$ and we work as above.

If $|A|\geq \o$, $I_0=\k \setminus A$
and $A=\bigcup _{n\in \N}I_n$ is any partition of $A$, then for each $F\in \Phi$ we have $F\cap I_n\neq \emptyset$, for all $n\in \N$.
\kdok
\begin{te}\label{T020}
{\rm ($\ast$)} If the size of the index set $I$ is less than the first measurable cardinal,
then the following conditions are equivalent
\begin{itemize}
\item[\rm (a)] The reduced product $\prod _\Phi \X _i$ is connected,
\item[\rm (b)] There are a finite set $K\subset I$ and $n\in \o$ such that
$\X _i$ is connected, for each $i\in K$,
and  $\{ i\in I: \X _i \models \f _n^{\rm conn}\}\cup K\in \Phi$.
\end{itemize}
The implication  (b) $\Rightarrow$ (a)  is true in ZFC (without any assumptions on $|I|$).
\end{te}
\dok
W.l.o.g.\ we can assume that $\rho _i \subset X_i$, $i\in I$, are reflexive relations.

(a) $\Rightarrow$ (b).
We prove the contrapositive: assuming  that
\begin{equation}\label{EQ051}
\forall K\in [I]^{<\o} \,\Big( \forall i\in K \;(\X_i\models \f ^{\rm conn} )
                               \Rightarrow \forall n\in \o \;\{ i\in I: \X _i \models \f _n^{\rm conn}\}\cup K\not\in \Phi \Big)
\end{equation}
we show that $\prod _\Phi \X _i$ is disconnected.

If the filter $\Phi$ contains a finite set,
then $\Phi$ is a principal filter generated by some $K\in [I]^{<\o}$
and, by (\ref{EQ051}), there is $i_0\in K$ such that $\X _{i _0}$ is disconnected.
Since $\Phi \upharpoonright K=\{ K\}$, by Fact \ref{T014} we have $\prod _\Phi \X _i\cong \prod _{\Phi \upharpoonright K} \X _i\cong \prod _{i\in K}\X _i$.
We take $x_{i_0},y_{i_0}\in X _{i_0}$ such that there is no path from $x_{i_0}$ to $y_{i_0}$ in $\X _{i_0}$
and choose arbitrary $x_i=y_i\in X_i$, for $i\in K\setminus \{ i_0\}$.
Assuming that there is a path from $x:=\la x_i:i\in K\ra$ to $y:=\la y_i:i\in K\ra$  in $\prod _{i\in K}\X _i$
we would obtain a path from $x_{i_0}$ to $y_{i_0}$ in $\X _{i_0}$, which is impossible.
Thus the reduced product $\prod _\Phi \X _i$ is disconnected.

So, the remaining case is when $\Phi$ is a filter consisting of infinite subsets of $I$.
Let $A_n :=\{ i\in I: \X _i \models \f ^{\rm conn}_n\}$, for $n\in \o$.
Then $A_0\subset A_1 \subset A_2 \subset \dots$
and, defining $I_0:=A_0$, $I_n:=A_{n}\setminus A_{n-1}$, for $n\in \N$, and $I_\infty=I\setminus \bigcup _{n\in \o}A_n$,
we have
$$\textstyle
I=\bigcup _{n\in\o}I _n \cup I_\infty .
$$
We continue with a case analysis.

1. $I_\infty \in \Phi$.
Then $|I_\infty|\geq \o$ is a non-measurable cardinal,
$\Phi\upharpoonright I_\infty $ is a filter on $I_\infty$ consisting of infinite sets,
and, by Fact \ref{T019},  there is a partition $I_\infty=\bigcup _{n\in \o}I_n'$ such that
for each $F\in \Phi\upharpoonright I_\infty$ we have $F\cap I_n'\neq \emptyset$, for infinitely many $n\in \o$.
If $n\in \o$, then for $i\in I_n'$ we have $i\not\in A_n$
and we choose $x_i,y_i\in X_i$ such that $d (x_i,y_i)>n+1$.
For $i\in I\setminus I_\infty$ we choose arbitrary $x_i=y_i\in X_i$.
Let $x:=\la x_i:i\in I\ra$, $y:=\la y_i:i\in I\ra$
and suppose that in $\prod _\Phi \X _i$ there is a path from $[x]_\Phi$ to $[y]_\Phi$
\begin{equation}\label{EQ052}
[x]_\Phi\r ^{\e_0}  [z^0]_\Phi\r ^{\e_1} [z^1]_\Phi \r ^{\e_2}
\dots \r ^{\e_{n-1}} [z^{n-1}]_\Phi\r ^{\e_n} [y]_\Phi.
\end{equation}
Then there would be $F\in \Phi$ which witnesses that; in other words,
\begin{equation}\label{EQ053}
\forall i\in F \;\;
x_i \r ^{\e_0}_i  z^0_i \r ^{\e_1}_i z^1_i \r ^{\e_2}_i \dots \r ^{\e_{n-1}}_i z^{n-1}_i\r ^{\e_n}_i y_i.
\end{equation}
Now let $m>n$ be such that there is $i\in F\cap I_m'$.
By (\ref{EQ053}) there is an ($n+1$)-path from $x_i$ to $y_i$, which is impossible by the choice of $x_i$ and $y_i$ ($d (x_i,y_i)>n+1$).
So, there is no path from $[x]_\Phi$ to $[y]_\Phi$ in $\prod _\Phi \X _i$
and this structure is disconnected.

2. $I_\infty \not\in \Phi$. Then we have $F\cap \bigcup _{n\in\o}I _n\neq\emptyset$, for all $F\in \Phi$.

2.1. Each $F\in \Phi$ intersects infinitely many sets $I_n$.
If $n\in \N$, then for $i\in I_n$ we have $i\not\in A_{n-1}$
and we choose $x_i,y_i\in X_i$ such that $d (x_i,y_i)>n$.
For $i\in I\setminus \bigcup _{n\in\N}I _n$ we choose arbitrary $x_i=y_i\in X_i$.
Let $x:=\la x_i:i\in I\ra$, $y:=\la y_i:i\in I\ra$
and suppose that in $\prod _\Phi \X _i$ there is a path from $[x]_\Phi$ to $[y]_\Phi$
\begin{equation}\label{EQ054}
[x]_\Phi\r ^{\e_0}  [z^0]_\Phi\r ^{\e_1} [z^1]_\Phi \r ^{\e_2}
\dots \r ^{\e_{n-1}} [z^{n-1}]_\Phi\r ^{\e_n} [y]_\Phi.
\end{equation}
Then there would be $F\in \Phi$ which witnesses that; in other words,
\begin{equation}\label{EQ055}
\forall i\in F \;\;
x_i \r ^{\e_0}_i  z^0_i \r ^{\e_1}_i z^1_i \r ^{\e_2}_i \dots \r ^{\e_{n-1}}_i z^{n-1}_i\r ^{\e_n}_i y_i.
\end{equation}
By the assumption of 2.1 there is $m>n+1$ such that $F\cap I_m\neq\emptyset$.
So, for $i\in F\cap I_m$, by (\ref{EQ055}) there is an ($n+1$)-path from $x_i$ to $y_i$,
which is impossible by the choice of $x_i$ and $y_i$ ($d (x_i,y_i)>m$).
So, there is no path from $[x]_\Phi$ to $[y]_\Phi$ in $\prod _\Phi \X _i$
and this structure is disconnected.

2.2. There are $n_0\in \o$ and $F_0\in \Phi$ such that $F_0\subset A_{n_0}\cup I_\infty$.
By (\ref{EQ051}) taking $K=\emptyset$ we obtain $A_n\not\in\Phi$, for all $n\in \o$.
So, assuming that $F\cap I_\infty=\emptyset$, for some $F\in \Phi$,
we would have $A_{n_0}\in\Phi$, which is false.
Thus
\begin{equation}\label{EQ056}
\forall F\in \Phi \;\;F\cap I_\infty\neq \emptyset,
\end{equation}
and $\Phi \upharpoonright I_\infty$ is a proper filter on $I_\infty $.

2.2.1. $|F\cap I_\infty|\geq \o$, for all $F\in \Phi$. Then $|I_\infty|\geq \o$ is a non-measurable cardinal,
$\Phi\upharpoonright I_\infty $ is a filter on $I_\infty$ consisting of infinite sets,
and, by Fact \ref{T019},  there is a partition $I_\infty=\bigcup _{n\in \o}I_n'$ such that
for each $F\in \Phi\upharpoonright I_\infty$ we have $F\cap I_n'\neq \emptyset$, for infinitely many $n\in \o$.
Now we continue exactly as in Case 1.

2.2.2. The set $F\cap I_\infty$ is finite, for some $F\in \Phi$.
Then $\Phi \upharpoonright I_\infty$ is a principal filter generated by some finite set $J\subset I_\infty$
and, hence, $A_{n_0}\cup J\in \Phi$.
So, by (\ref{EQ051}), there is $i_0\in J$ such that $\X _{i_0}$ is disconnected.
We take $x_{i_0},y_{i_0}\in X _{i_0}$ such that there is no path from $x_{i_0}$ to $y_{i_0}$ in $\X _{i_0}$
and choose arbitrary $x_i=y_i\in X_i$, for $i\in I\setminus \{ i_0\}$.
Let $x:=\la x_i:i\in I\ra$, $y:=\la y_i:i\in I\ra$
and suppose that there is a path from $[x]_\Phi$ to $[y]_\Phi$ in $\prod _\Phi \X _i$,
\begin{equation}\label{EQ402}
[x]_\Phi\r ^{\e_0}  [z^0]_\Phi\r ^{\e_1} [z^1]_\Phi \r ^{\e_2}
\dots \r ^{\e_{n-1}} [z^{n-1}]_\Phi\r ^{\e_n} [y]_\Phi.
\end{equation}
Then there would be $F\in \Phi$ which witnesses that; in other words,
\begin{equation}\label{EQ058}
\forall i\in F \;\;
x_i \r ^{\e_0}_i  z^0_i \r ^{\e_1}_i z^1_i \r ^{\e_2}_i \dots \r ^{\e_{n-1}}_i z^{n-1}_i\r ^{\e_n}_i y_i.
\end{equation}
But $i_0\in J \subset F$
and we obtain a path from $x_{i_0}$ to $y_{i_0}$ in $\X _{i_0}$, which is impossible.
Thus the reduced product $\prod _\Phi \X _i$ is disconnected.

(b) $\Rightarrow$ (a). Let $K\in [I]^{<\o}$, $n\in \o$ and, defining $A:=\{ i\in I: \X _i \models \f _n^{\rm conn}\}$, suppose that

(i) $\X _i$ is connected, for each $i\in K$, and

(ii) $A \cup K\in \Phi$.

\noindent
Since $\f _n^{\rm conn}\Rightarrow \f ^{\rm conn}$ we can assume that $A\cap K=\emptyset$.
For $x=\la x_i :i\in I\ra, y=\la y_i :i\in I\ra \in \prod _{i\in I}X_i$
we construct a path  from $[x]_\Phi$ to $[y]_\Phi$ in $\prod _{\Phi }\X _i$.

If $i\in A$, then $\X _i \models \f _n^{\rm conn}$
and, since $x_i,y_i\in X_i$,
there is a path $\la \bar z,\e\ra$ of length $n+1$ from $x_i$ to $y_i$;
thus we can choose one such $\e$.
Let $\e ^k \in 2^{n+1}$, $k<l$, be an enumeration of the $n+1$-tuples
chosen in this way for $i\in A$.
Then defining $A_{\e ^k}=\{ i\in A : \e ^k \mbox{ is chosen for } x_i \mbox{ and } y_i \}$, for $k<l$, we have
$$\textstyle
A=\bigcup _{k<l} A_{\e ^k} .
$$
Let $y'\in \prod _{i\in I}X_i$ be defined by
\begin{equation}\label{EQ047}
y'_i=\left\{\begin{array}{ll}
                   x_i & \mbox{ if } i\in I\setminus A ,\\
                   y_i & \mbox{ if } i\in A.
               \end{array}
         \right.
\end{equation}
We define sequences $s^m=\la s^m_i :i\in I\ra \in \prod _{i\in I}X_i$, for $m\leq l-2$,
such that there are $l$-many paths of length $(n+1)$ witnessing that
$[x]_\Phi \sim [s^0]_\Phi \sim [s^1]_\Phi\sim \dots \sim [s^{l-2}]_\Phi \sim [y']_\Phi$, in the following way
\begin{equation}\label{EQ025}
s^m_i=\left\{\begin{array}{ll}
                   x_i & \mbox{ if } i\in I\setminus \bigcup _{k\leq m}A_{\e ^k} ,\\
                   y_i & \mbox{ if } i\in \bigcup _{k\leq m}A_{\e ^k}.
               \end{array}
         \right.
\end{equation}
Now, $x$ and $s^0$ satisfy the assumptions of Fact \ref{T017} over $A_{\e ^0}$;
so, $[x]_\Phi \sim [s^0]_\Phi$.
$s^0$ and $s^1$ satisfy the assumptions of Fact \ref{T017} over $A_{\e ^1}$;
so, $[s^0]_\Phi \sim [s^1]_\Phi$.
Continuing in this way we have that
$s^{l-2}$ and $y'$ satisfy the assumptions of Fact \ref{T017} over $A_{\e ^{l-1}}$;
so, $[s^{l-2}]_\Phi \sim [y']_\Phi$. Thus we have  $[x]_\Phi \sim [y']_\Phi$.

Let $K=\{ i_j:j<k\}$ be an enumeration.
For $j<k$ the structure $\X _{i_j}$ is connected
so we can choose a path $\la \bar z _j,\e _j\ra$ of length $n_j$ from $x_{i_j}$ to $y_{i_j}$.
By (ii) we have $F:=A\cup K \in \Phi$.
Let $y''\in \prod _{i\in I}X_i$ be defined by
\begin{equation}\label{EQ403}
y''_i=\left\{\begin{array}{ll}
                   x_i & \mbox{ if } i\in I\setminus F ,\\
                   y_i & \mbox{ if } i\in F.
               \end{array}
         \right.
\end{equation}
We define sequences $t^j=\la t^j_i :i\in I\ra \in \prod _{i\in I}X_i$, for $j\leq k-2$,
such that $[y']_\Phi \sim [t^0]_\Phi \sim [t^1]_\Phi\sim \dots \sim [t^{k-2}]_\Phi \sim [y'']_\Phi$, by
\begin{equation}\label{EQ404}
t^j_i=\left\{\begin{array}{ll}
                   x_i & \mbox{ if } i\in I\setminus (A\cup \{ i_0 ,\dots, i_j\}),\\
                   y_i & \mbox{ if } i\in A\cup \{ i_0 ,\dots, i_j\}.
               \end{array}
         \right.
\end{equation}
Namely, $y'$ and $t^0$ satisfy the assumptions of Fact \ref{T017} over $\{i_0\}$;
so, $[y']_\Phi \sim [t^0]_\Phi$.
$t^0$ and $t^1$ satisfy the assumptions of Fact \ref{T017} over $\{i_1\}$;
so, $[t^0]_\Phi \sim [t^1]_\Phi$.
Continuing in this way we have that
$t^{k-2}$ and $y''$ satisfy the assumptions of Fact \ref{T017} over $\{ i_{k-1}\}$;
so $[t^{k-2}]_\Phi \sim [y'']_\Phi$.
Thus we have proved that $[y']_\Phi \sim  [y'']_\Phi$
and, since $[x]_\Phi \sim [y']_\Phi$ and $[y'']_\Phi =[y]_\Phi$,
we finally have $[x]_\Phi \sim [y]_\Phi$.
\kdok
\begin{rem} \rm
Regarding Theorem \ref{T020}, if $|I|\geq \o$ and the filter $\Phi$ contains all cofinite sets, then condition (b) can be replaced by the following simpler condition
\begin{itemize}
\item[\rm (b')] There is $n\in \o$ such that  $\{ i\in I: \X _i \models \f _n^{\rm conn}\}\in \Phi$.
\end{itemize}
\end{rem}

In particular cases, when the filter $\Phi$ is minimal (i.e.\ $\Phi =\{ I\}$) and when the filter $\Phi$ is maximal (an ultrafilter in $P(I)$)
by Theorem \ref{T020} we have the following two corollaries.
\begin{cor}\label{T401}
(Direct products and powers)
If {\rm ($\ast$)} holds and  $\Phi =\{ I\}$, then

(a) $\prod _{i\in I} \X _i$ is connected
iff all the structures $\X _i$ are connected and there is $n\in \o$ such that  $\X _i $ is of diameter $\leq n$, for all except finitely many $i\in I$;

(b) The direct power $\prod _{i\in I} \X $ is connected iff $\X$ is of finite diameter;

(c) If $|I|<\o$, then $\prod _{i\in I} \X _i$ is connected iff all the structures $\X _i$ are connected.
\end{cor}
\begin{cor}\label{T402}
(Ultraproducts and ultrapowers)
If {\rm ($\ast$)} holds, $\o \leq |I|$ is less than the first measurable cardinal and $\CU \subset P(I)$ is a non-principal ultrafilter, then

(a) $\prod _\CU \X _i$ is connected, iff $\{ i\in I: \X _i \models \f _n^{\rm conn}\}\in \CU$, for some $n\in \o$;

(b) The ultrapower $\prod _\CU \X $ is connected iff $\X$ is of finite diameter.

\noindent
The implications  $\Leftarrow$  are true in ZFC (without any assumptions on $|I|$).
\end{cor}
\begin{te}\label{T022}
Let $|I|$ be a measurable cardinal, $\CU$ a non-principal $\o$-complete ultrafilter on $I$ and $\X _i=\la X_i ,\r _i\ra$, $i\in I$, $L$-structures. Then
$$\textstyle
\prod _\CU \X _i \mbox{ is connected} \Leftrightarrow \{ i\in I : \X _i \mbox{ is connected}\}\in \CU .
$$
\end{te}
\dok
($\Leftarrow$)
Let $U:=\{ i\in I : \X _i \mbox{ is connected}\}\in \CU$ and let
$x=\la x_i :i\in I\ra, y=\la y_i :i\in I\ra \in \prod _{i\in I}X_i$.
We construct a path  from $[x]_\CU$ to $[y]_\CU$ in $\prod _{\CU }\X _i$.
For $i\in U$ there is a path $\la \bar z ^i,\e ^i\ra$ of length $n_i+1$ from $x_i$ to $y_i$
and we  choose such $\e _i \in 2^{n_i+1}$.
Since $|\bigcup _{n\in \o}2^{n+1}|=\o$, there are $\leq \o$ such choices;
let $\{ \e^k:k<\lambda \}$ be an enumeration of them; so, $\lambda \leq \o$.
Then we have $U=\bigcup _{k< \lambda}A_k$, where $A_k:=\{i\in U : \e ^i=\e ^k\}$
and, since $\CU$ is $\o$-complete, there is $k<\lambda$ such that $A_k \in \CU$.
Now, using Fact \ref{T017} we construct a path  from $[x]_\CU$ to $[y]_\CU$.

($\Rightarrow$)
Let $U:=\{ i\in I : \X _i \mbox{ is not  connected}\}\in \CU$.
For each $i\in U$ we pick $x_i,y_i\in X_i$ such that $d(x_i,y_i)=\infty$,
for $i\in I\setminus U$ we take arbitrary $x_i=y_i\in X_i$
and define $x=\la x_i :i\in I\ra, y=\la y_i :i\in I\ra \in \prod _{i\in I}X_i$.
Assuming that there is a path   from $[x]_\CU$ to $[y]_\CU$ in $\prod _\CU \X _i $,
\begin{equation}\label{EQ022}
[x]_\CU\r ^{\e_0}  [z^0]_\CU\r ^{\e_1} [z^1]_\CU \r ^{\e_2}
\dots \r ^{\e_{n-1}} [z^{n-1}]_\CU \r ^{\e_n} [y]_\CU,
\end{equation}
there would be $U_1\in \CU$ which witnesses that and, hence,
\begin{equation}\label{EQ023}
\forall i\in U\cap U_1 \;\;
x_i \r ^{\e_0}_i  z^0_i \r ^{\e_1}_i z^1_i \r ^{\e_2}_i \dots \r ^{\e_{n-1}}_i z^{n-1}_i\r ^{\e_n}_i y_i.
\end{equation}
which is impossible because $d (x_i,y_i)=\infty$, for all $i\in U$.
So, there is no path from $[x]_\CU$ to $[y]_\CU$ in $\prod _\CU \X _i$
and this structure is disconnected.
\kdok
A simple example of a connected $L$-structure of infinite diameter
is the linear graph $\G_\o =\la \o ,\r\ra$, where $\r =\{ \la m,n \ra : |m-n|=1\}$.
In the absence of measurable cardinals Theorem \ref{T020} gives a characterization of connected reduced products of $L$-structures.
In the presence of measurable cardinals we have the following equivalence.
\begin{te}\label{T400}
For each infinite set $I$ the following conditions are equivalent

(i) The cardinal $|I|$ is less than the first measurable cardinal,

(ii) The implication (a) $\Rightarrow$ (b) of Theorem \ref{T020} is true for each filter $\Phi \subset P(I)$ and each family of $L$-structures $\{ \X_i :i\in I\}$,

(iii) The ultraproduct $\prod _{\CU}\G_\o$ is a disconnected graph for each non-principal ultrafilter $\CU \subset P(I)$.
\end{te}
\dok
The implication (i) $\Rightarrow$ (ii) is proved in Theorem \ref{T020}.

(ii) $\Rightarrow$ (iii). Let (ii) hold and let $\CU \subset P(I)$ be a non-principal ultrafilter. Assuming that the ultraproduct $\prod _{\CU}\G_\o$ is connected
we would have a finite set $K\subset I$ and $n\in \o$ such that $\{ i\in I: \G _\o \models \f _n^{\rm conn}\}\cup K =K\in \Phi$, which is false.

(iii) $\Rightarrow$ (i). Let (iii) hold.
Suppose that $|I|\geq \k$, where $\k$ is the first measurable cardinal
and w.l.o.g.\ suppose that $\k\subset I$.
Let $\CU _\k\subset P(\k )$ be a non-principal $\o$-complete ultrafilter
and let $\CU \subset P(I)$ be its extension to $I$.
Then by Fact \ref{T014} we have $\prod _{\CU}\G_\o :=(\prod _{i\in I}\G_\o)/\CU \cong (\prod _{i\in \k }\G_\o)/\CU _\k$
and by Theorem \ref{T022} the ultraproduct $(\prod _{i\in \k }\G_\o)/\CU _\k$ is connected,
which is false by (iii).
\kdok
The connectivity relation (and, hence, the components) of a reduced product are described by the following theorem.
\begin{te}\label{T403}
{\rm ($\ast$)} If  $x=\la x_i:i\in I\ra, y=\la y_i:i\in I\ra \in \prod _{i\in I}X_i$, then
$$
[x]_\Phi \sim [y]_\Phi \;\;\Leftrightarrow\;\; \exists n\in \o \;\;\{ i\in I: \X _i \models \f _{d\leq n+1}[x_i,y_i]\}\in \Phi.
$$
\end{te}
\dok
($\Rightarrow$) If $[x]_\Phi \sim [y]_\Phi$, then in $\prod _\Phi \X _i$ there is a path from $[x]_\Phi$ to $[y]_\Phi$;
namely, there are $n\in \o$, $\e \in 2^{n+1}$ and $z^0,\dots,z^{n-1}\in \prod _{i\in I}X_i$ such that we have
$[x]_\Phi\r ^{\e_0}  [z^0]_\Phi\r ^{\e_1} [z^1]_\Phi \r ^{\e_2}\dots \r ^{\e_{n-1}} [z^{n-1}]_\Phi\r ^{\e_n} [y]_\Phi$.
Since $\Phi$ is closed under finite intersections there is $F\in \Phi$ such that
for each $i\in F$ we have $x_i \r ^{\e_0}_i  z^0_i \r ^{\e_1}_i z^1_i \r ^{\e_2}_i \dots \r ^{\e_{n-1}}_i z^{n-1}_i\r ^{\e_n}_i y_i$.
So $F\subset \{ i\in I: \X _i \models \f _{d\leq n+1}[x_i,y_i]\}\in \Phi$.

($\Leftarrow$) Let $n\in \o$ and $A:=\{ i\in I: \X _i \models \f _{d\leq n+1}[x_i,y_i]\}\in \Phi$.
If $i\in A$, then there is a path $\la \bar z,\e\ra$ of length $n+1$ from $x_i$ to $y_i$;
thus we can choose one such $\e$.
Let $\e ^k \in 2^{n+1}$, $k<l$, be an enumeration of the $n+1$-tuples
chosen in this way for $i\in A$.
Then defining $A_{\e ^k}=\{ i\in A : \e ^k \mbox{ is chosen for } x_i \mbox{ and } y_i \}$, for $k<l$,
we have $A=\bigcup _{k<l} A_{\e ^k}$.
Let $y'\in \prod _{i\in I}X_i$, where $y'\upharpoonright A=y\upharpoonright A$ and $y'\upharpoonright (I\setminus  A)=x\upharpoonright (I\setminus  A)$.
Then $[y']_\Phi =[y]_\Phi$. We define sequences $s^m=\la s^m_i :i\in I\ra \in \prod _{i\in I}X_i$, for $m\leq l-2$,
such that there are $l$-many paths of length $(n+1)$ witnessing that
$[x]_\Phi \sim [s^0]_\Phi \sim [s^1]_\Phi\sim \dots \sim [s^{l-2}]_\Phi \sim [y']_\Phi$, by
\begin{equation}\label{EQ401}
s^m_i=\left\{\begin{array}{ll}
                   x_i & \mbox{ if } i\in I\setminus \bigcup _{k\leq m}A_{\e ^k} ,\\
                   y_i & \mbox{ if } i\in \bigcup _{k\leq m}A_{\e ^k}.
               \end{array}
         \right.
\end{equation}
Now, $x$ and $s^0$ satisfy the assumptions of Fact \ref{T017} over $A_{\e ^0}$;
so, $[x]_\Phi \sim [s^0]_\Phi$.
$s^0$ and $s^1$ satisfy the assumptions of Fact \ref{T017} over $A_{\e ^1}$;
so, $[s^0]_\Phi \sim [s^1]_\Phi$.
Continuing in this way we have that
$s^{l-2}$ and $y'$ satisfy the assumptions of Fact \ref{T017} over $A_{\e ^{l-1}}$;
so, $[s^{l-2}]_\Phi \sim [y']_\Phi$. Thus we have  $[x]_\Phi \sim [y']_\Phi=[y]_\Phi$.
\kdok

\paragraph{Acknowledgement} This research was supported by the Science Fund of the Republic of Serbia, Program PROMIS, Grant No.\ 6062228, {\it Classification of Large Objects--Ultrafilters and Directed Sets}--CLOUDS.

\end{document}